\documentclass{ifacconf}

\usepackage{amsmath}
\usepackage{amssymb}
\usepackage{latexsym}
\usepackage{amsfonts}
\usepackage{graphicx}      
\usepackage{xcolor}
\usepackage{natbib}        

\theoremstyle{definition}
\newtheorem{ex}[thm]{Example}

\newcommand{\R}{\ensuremath{\mathbb R}}    
\newcommand{\N}{\ensuremath{\mathbb N}}    

\newcommand{\<}{\langle}
\renewcommand{\>}{\rangle}

\newcommand{\calH}{\mathcal H}

\newcommand{\calR}{\mathcal R}

\newcommand{\calU}{\mathcal U}

\renewcommand{\Re}{\operatorname{Re}}

\newcommand{\dom}{\operatorname{dom}}

\newcommand{\Sra}{\Rightarrow}

\newcommand{\Sla}{\Leftarrow}

\newcommand{\ol}{\overline}

\newcommand{\wt}{\widetilde}

\newcommand{\tf}{\mathbf{t}}

\newcommand{\ph}{port-Hamitonian }

\begin{document}
\begin{frontmatter}

\title{Representing the dissipation of infinite-dimensional linear port-Hamiltonian systems} 

\thanks[footnoteinfo]{FP was funded by the Carl Zeiss Foundation within the project {\it DeepTurb--Deep Learning in and from Turbulence}. He was further supported by the free state of Thuringia and the German Federal Ministry of Education and Research (BMBF) within the project {\it THInKI--Th\"uringer Hochschulinitiative für KI im Studium}.}

\author[First]{Friedrich M.\ Philipp}

\address[First]{Technische Universität Ilmenau, Optimization-based Control Group, Institute of Mathematics, 98693 Ilmenau, Germany (e-mail: friedrich.philipp@tu-ilmenau.de).}

\begin{abstract}                
It is well known that linear and non-linear dissipative port-Hamiltonian systems in finite dimensions admit an energy balance, relating the energy increase in the system with the supplied energy and the dissipated energy. The integrand in the dissipation term is then a function of the state variable. In this note, we answer the question of when this is possible for linear port-Hamiltonian systems in infinite dimensions.
\end{abstract}


\end{frontmatter}

\section{Introduction}
The port-Hamiltonian modeling paradigm is so appealing as it provides an elegant mathematical framework for energy-based modeling and analysis of multi-physics systems, see, e.g., \cite{SchaJelt2014}. Port-Hamiltonian systems, providing in addition various coupling capabilities, are therefore very well-suited to describe the energy flows, energy conservation and interconnection of physical systems in a wide range of applications, such as robotics (\cite{AngeMusi17}), renewable energy systems (\cite{TonsKapa23}), and adaptive building (\cite{SchaZell24}), just to name a few.

A conventional input-state-output \ph ODE-system has the form
\begin{align}
\begin{split}
\dot x &= A(x)\nabla H(x) + B(x)u\\
y &= B(x)^*\nabla H(x),
\end{split}
\end{align}
where $H$ is the system's Hamiltonian describing the energy in the system and $A$ is pointwise dissipative, i.e., $\Re A(x) := \frac 12(A(x)+A(x)^*)$ is negative semi-definite. Solutions of the system can be easily shown to satisfy the energy balance relation, which is given by
\begin{align*}
H(x(t_2)) &- H(x(t_1))\\
&= \int_{t_1}^{t_2}y^*u\,dt - \int_{t_1}^{t_2}\|[-\Re A(x)]^{1/2}\nabla H(x)\|^2\,dt,
\end{align*}
where $t_2>t_1\ge 0$. While the left-hand side describes the increase of energy in the system on the time horizon $[t_1,t_2]$, on the right-hand side we have the difference between the energy supplied to the system $\int_{t_1}^{t_2}y^*u\,dt$ and the dissipated energy $\int_{t_1}^{t_2}\|[\Re A(x)]^{1/2}\nabla H(x)\|^2\,dt$. In the special linear case, where $H(x) = \frac 12\|x\|_2^2$ and $A(x)\equiv A$, the dissipation term has the simple form $\int_{t_1}^{t_2}\|R^{1/2}x(t)\|^2\,dt$, where $R := -\Re A = -\frac 12(A+A^*)\ge 0$.

If the state space is an infinite-dimensional Hilbert space and $A$ is a bounded operator, the above definition of the dissipation operator $R$ still makes sense. However, in applications the operator $A$ is usually unbounded, and $A+A^*$ is defined on the intersection of subspaces $\dom A\cap\dom A^*$. This subspace is the object of investigation in \cite{ArliTret20}, where it is proved in particular that there even exist generators $A$ of analytic $C_0$-semigroups for which this intersection is trivial: $\dom A\cap\dom A^* = \{0\}$. One of the results of the present note is that for such generators there however exists a self-adjoint operator $R$ such that the energy balance holds for mild solutions $x$ with the dissipation term $\int_{t_1}^{t_2}\|R^{1/2}x(t)\|^2\,dt$.

Infinite-dimensional linear port-Hamiltonian systems have been investigated in a variety of works (see, e.g., \cite{Augner19,Jacob2019,Jasc23,RashCali20,Vill07}). The standard book on the topic is surely \cite{Jacob2012}, which mainly covers boundary control problems of port-Hamiltonian PDEs on a one-dimensional spacial domain. This setting has been generalized to multi-dimensional domains in \cite{Skre21}. Recently, a wide class of linear infinite-dimensional port-Hamiltonian systems has been introduced in \cite{PhilReis23}, which is based on the framework of Staffan's system nodes (\cite{Staf05}) and admits both boundary and distributed control. For the sake of simplicity, we restrict ourselves to distributed control, here.

The paper is organized as follows. In Section \ref{s:ph}, we briefly present the setting of the paper and discuss the problem of representing the dissipation term in the energy balance relation. Next, in Section \ref{s:icc}, it is shown that there always exists a linear operator $\calR$, mapping initial value $x_0$ and control $u$ to a function $\calR(x_0,u)(t)$ such that the dissipation rate can be represented as $\|\calR(x_0,u)(t)\|_A^2$, where $\|\cdot\|_A$ denotes the graph norm of $A$. Finally, in Section \ref{s:full} we prove that there exists a self-adjoint operator $R$ such that the dissipation rate can be written in the form $\|R^{1/2}x\|^2$ if and only if the quadratic form $x\mapsto\Re\<Ax,x\>$ is closable.

\section{Linear input-state-output port-Hamiltonian systems}\label{s:ph}
In what follows, let $\calH$ and $\calU$ be a Hilbert spaces. We shall consider infinite-dimensional linear control systems of the following form:
\begin{align}
\dot x &= Ax + Bu,\quad x(0) = x_0,\label{e:dynamics}\\
y &= B^*x.\label{e:output}
\end{align}
These will be termed ({\em input-state-output}) {\em port-Hamiltonian systems}. Hereby, we assume that
\begin{itemize}
	\item $A$ is a maximal dissipative\footnote{i.e., $A$ is the generator of a $C_0$-semigroup of contractions. By the Lumer-Phillips theorem, this is the case if and only if $\Re\<Ax,x\>\le 0$ for all $x\in\dom A$ (i.e., $A$ is {\em dissipative}) and $A-I$ is surjective.} operator in $\calH$,
	\item $B\in L(\calU,\calH)$.
\end{itemize}
If $u\in L^1_{\rm loc}(\R_0^+,\calU)$ and $x_0\in\calH$, the mild solution $x_u(\,\cdot\,;x_0)$ of \eqref{e:dynamics} is defined by
$$
x_u(t;x_0) = S(t)x_0 + \int_0^t S(t-s)Bu(s)\,ds,
$$
where $(S(t))_{t\ge 0}$ denotes the $C_0$-semigroup of contractions generated by $A$.

If $u\in C^1(\R_0^+,\calU)$ and $x_0\in\dom A$, the mild solution $x_u(\,\cdot\,;x_0)$ will be called the classical solution of \eqref{e:dynamics}. In this case, we have $x\in C^1(\R_0^+,\calH)$ and $x(t)\in\dom A$ for all $t\ge 0$, cf.\ Corollaries VI.7.6 and VI.7.8 in \cite{EngelNagel}.

The {\em energy Hamiltonian} $H : \calH\to\R$ of the system \eqref{e:dynamics}--\eqref{e:output} is defined by
$$
H(x) := \tfrac 12\|x\|^2.
$$
If $x$ is a classical solution of \eqref{e:dynamics}, then
\begin{align*}
\tfrac d{dt}H(x)
&= \<\dot x,x\> = \Re\<Ax+Bu,x\>\\
&= \Re\<u,y\> + \Re\<Ax,x\>,
\end{align*}
and hence, for $T>0$,
\begin{equation}\label{e:eb_classical}
H(x(T)) - H(x(0)) = \int_{0}^{T}\Re\<u,y\>\,dt + \int_{0}^{T}\Re\<Ax,x\>\,dt,
\end{equation}
which is the well-known {\em energy balance}. Moreover, if $\calH$ is finite-dimensional (and hence $A$ is a dissipative matrix) we may write $\Re\<Ax,x\> = -\<Rx,x\> = -\|R^{1/2}x\|^2$, where $R = -\Re A = -\frac 12(A+A^*)$ is the real part of $-A$. In the infinite-dimensional case, we face two problems:
\begin{itemize}
	\item The energy balance \eqref{e:eb_classical} holds for classical solutions $x$ of \eqref{e:dynamics}; however, for mild solutions $x$, it does not make sense because $x(t)\notin\dom A$, in general;\\[-5pt]
	\item Attempting to write the dissipation rate $-\Re\<Ax,x\>$ as $\|R^{1/2}x\|^2$ with the real part $R$ of $-A$ fails in general since $R = -\frac 12(A+A^*)$ is symmetric, but might not be self-adjoint. Its domain of definition $\dom R = \dom A\cap\dom A^*$ might even be the trivial subspace $\{0\}$. On the other hand, Example \ref{ex:ex} exposes a maximal dissipative operator $A$ with $\dom R$ being dense in $\calH$, but $R=0$ on its domain.
\end{itemize}
In this note it is our aim to replace the expression $\Re\<Ax,x\>$ by other ones which are also meaningful for mild solutions $x$. To this end, we introduce the symmetric sesquilinear form $\mathbf r$ with domain $\dom(\mathbf r) = \dom A$, defined by
$$
\mathbf r[x,y] := -\tfrac 12\big(\<Ax,y\> + \<x,Ay\>\big),\qquad x,y\in\dom A.
$$
The form $\mathbf r$ is called the {\em real part} of the form $(x,y)\mapsto\<Ax,y\>$. Note that $\mathbf r[x] := \mathbf r[x,x] = -\Re\<Ax,x\>\ge 0$ for $x\in\dom A$.

\section{Representation in terms of initial condition and control}\label{s:icc}
In what follows, let $\|x\|_A := (\|x\|^2 + \|Ax\|^2)^{1/2}$ for $x\in\dom A$ and
$$
\calH_A = (\dom A,\|\cdot\|_A).
$$

\smallskip
\begin{prop}\label{p:timo}
There exists a linear operator
$$
\calR : \calH\times L^2_{\rm loc}(\R_0^+,\calU)\to L^2_{\rm loc}(\R_0^+,\calH_A)
$$
such that for $x_0\in\calH$, $u\in L^2_{\rm loc}(\R_0^+;\calU)$, the associated mild solution $x(t) = x_u(t;x_0)$, and $T>0$ we have
\begin{equation}\label{e:eb_general}
H(x(T)) - H(x_0) = \int_0^T\Re\<u,y\>\,dt - \int_0^T\|\calR(x_0,u)(t)\|_{A}^2\,dt.
\end{equation}
If $u_1=u_2$ on $[0,T]$, then also $\calR(x_0,u_1)=\calR(x_0,u_2)$ on $[0,T]$ for all $x_0\in\calH$, and the restriction operator $\calR_T : \calH\times L^2(0,T;\calU)\to L^2(0,T;\calH_A)$ is bounded.
\end{prop}
\begin{pf}
For $x\in\dom A$ we have $\mathbf r[x]\le\|Ax\|\|x\|\le\frac 12\|x\|_A^2$, $x\in\dom A$. Hence, there exists a bounded non-negative self-adjoint operator $M$ in $\calH_A$ such that $\mathbf r[x,y] = \<Mx,y\>_A$ for $x,y\in\dom A$.

Let $T>0$ be fixed but arbitrary. For $(x_0,u)\in\dom A\times C^1(0,T;\calU)$ define
$$
\calR_T(x_0,u)(t) := M^{1/2}x_u(t;x_0),\quad t\in [0,T].
$$
Let $x(t) := x_u(t;x_0)$. Then $\|\calR_T(x_0,u)(t)\|_{A}^2 = \mathbf r[x(t)]$ and thus
\begin{align*}
\int_0^T\|&\calR_T(x_0,u)(t)\|_A^2\,dt\\
&= \int_0^T\Re\<u,y\>\,dt - \tfrac 12\|x(T)\|^2 + \tfrac 12\|x_0\|^2\\
&\le \int_0^T\Re\<Bu(t),x(t)\>\,dt + \tfrac 12\|x_0\|^2\\
&\le \|B\|\int_0^T\|u(t)\|\|x(t)\|\,dt + \tfrac 12\|x_0\|^2\\
&\le \|B\|\|u\|_{L^2(0,T;\,\calU)}\|x\|_{L^2(0,T;\calH)} + \tfrac 12\|x_0\|^2.
\end{align*}
Now, as
\begin{align*}
\|x(t)\|
&\le \|S(t)x_0\| + \int_0^t \|S(t-s)Bu(s)\|\,ds\\
&\le \|x_0\| + \|B\|\int_0^t\|u(s)\|\,ds,
\end{align*}
we obtain
$$
\|x\|_{L^2(0,T;\calH)}\,\le\,\sqrt{2T}\|x_0\| + T\|B\|\|u\|_{L^2(0,T;\calU)}.
$$
Hence, one easily sees that
\begin{align*}
\|\calR_T(x_0,u)\|_{L^2(0,T;\calH_A)}
&\le \sqrt{T}\|B\|\|u\|_{L^2(0,T;\,\calU)} + \frac 1{\sqrt 2}\|x_0\|
\end{align*}
for all $(x_0,u)\in\dom A\times C^1(0,T;\calU)$. Therefore, and since $\dom A\times C^1(0,T;\,\calU)$ is dense in $\calH\times L^2(0,T;\,\calU)$, the linear operator $R_T$ uniquely extends to a bounded operator on $\calH\times L^2(0,T;\,\calU)$, which we also denote by $\calR_T$. The identity \eqref{e:eb_general} then extends by boundedness to $(x_0,u)\in \calH\times L^2(0,T;\,\calU)$.

Finally, it is easily seen that for $T_2 > T_1 > 0$ and $u\in L^2(0,T_2;\calU)$ we have
$$
\calR_{T_2}(x_0,u)|_{[0,T_1]} = \calR_{T_1}(x_0,u|_{[0,T_1]}).
$$
Therefore, for $u\in L^2_{\rm loc}(\R_0^+,\calU)$ it makes sense to define $\calR(x_0,u)(t) := \calR_T(x_0,u|_{[0,T]})(t)$ for any $T\ge t$.\hfill $\square$
\end{pf}


\medskip
\begin{ex}\label{ex:ex}
Let $Ax = x'$, $\dom A = \{x\in H^1(0,1) : x(1)=0\}$, in $\calH = L^2(0,1)$. Then $A$ generates the contractive $C_0$-semigroup $S(\cdot)$, given by
\[
[S(t)x](\omega) = 
\begin{cases}
x(t+\omega) &\text{if $t+\omega\le 1$}\\
0 &\text{otherwise}.
\end{cases}
\]
It is easily verified that the operator $M^{1/2}$ from the proof of Proposition \ref{p:timo} is given by
\[
M^{1/2}x = \frac{e\sqrt 2}{\sqrt{e^4-1}}x(0)\sinh(1-\,\cdot\,),\quad x\in \dom A.
\]
Hence,
\begin{align*}
\calR(x_0,u)(t) = \frac{e\sqrt 2}{\sqrt{e^4-1}}x_u(t;x_0)(0)\cdot\sinh(1-\,\cdot\,),
\end{align*}
where (with $m(t) := \max\{t-1,0\}$)
\[
x_u(t;x_0)(0) = \chi_{[0,1]}(t)x_0(t) + \int_{m(t)}^t [Bu(s)](t-s)\,ds.
\]
Since $\|\sinh(1-\,\cdot\,)\|_A = \frac 12\frac{\sqrt{e^4-1}}{e}$, for the dissipation rate we get $\|\calR(x_0,u)(t)\|_A^2 = \frac 12|x_u(t;x_0)(0)|^2$. And indeed,
\[
-\Re\<Ax_u(t;x_0),x_u(t;x_0)\> = \tfrac 12|x_u(t;x_0)(0)|^2.
\]
In particular, this shows that in this special case the operator $\calR$ can be replaced by $\wt\calR$, mapping to $L^2(\R_0^+)$, where $\wt\calR(x_0,u) = 2^{-1/2}\cdot x_u(\,\cdot\,;x_0)(0)$. Furthermore, note that $A+A^* = 0$ on its domain $\{x\in H^1(0,1) : x(0)=x(1)=0\}$ so that this operator does not provide any information on the dissipation.
\end{ex}

\section{Representation in terms of the state}\label{s:full}
In the finite-dimensional case, we have that $-\Re\<Ax,x\> = \|R^{1/2}x\|^2$, where $R = \frac 12(A+A^*)$. In this section we consider the case where there exists a non-negative (possibly unbounded) self-adjoint operator $R$ in $\calH$ with $\dom R^{1/2}\supset\dom A$ such that for $x_0\in\calH$ and $u\in L^2_{\rm loc}(\R_0^+,\calU)$ we have
\begin{equation}\label{e:RR}
\|\calR(x_0,u)(t)\|_A = \|R^{1/2}x_u(t;x_0)\|\qquad\text{for a.e.\ $t\ge 0$}.
\end{equation}
Recall that a non-negative symmetric form $\mathbf a$ on $\calH$ is said to be {\em closed}, if for any sequence $(x_n)\subset\dom(\mathbf a)$ with $x_n\to x$ ($n\to\infty$) and $\mathbf a[x_n-x_m]\to 0$ ($n,m\to\infty$) we have $x\in\dom(\mathbf a)$ and $\mathbf a[x_n]\to \mathbf a[x]$ ($n\to\infty$). The form $\mathbf a$ is called {\em closable} if it admits a closed extension. Equivalently, $\mathbf a$ is closable if $x_n\to 0$ and $\mathbf a[x_n-x_m]\to 0$ as $n,m\to\infty$ imply $\mathbf a[x_n]\to 0$ as $n\to\infty$.

\begin{thm}
There exists a non-negative self-adjoint operator $R$ in $\calH$ such that \eqref{e:RR} holds for all $x_0\in\calH$ and $u\in L^2_{\rm loc}(\R_0^+,\calU)$ if and only if the form $\mathbf r$ is closable.
\end{thm}
\begin{pf}
``$\Sra$'' Let the form $\mathbf r$ be closable. Denote by $\ol{\mathbf r}$ the closure of $\mathbf r$. Then, by the second representation theorem (see \cite{kato}), there exists a non-negative self-adjoint operator $R$ in $\calH$ such that $\dom(\ol{\mathbf r}) = \dom R^{1/2}$ and
$$
\ol{\mathbf r}[x,y] = \<R^{1/2}x,R^{1/2}y\>\quad\text{for all $x,y\in\dom(\ol{\mathbf r})$}.
$$
In particular, we have $\dom A\subset\dom R^{1/2}$. Let $u\in C^1([0,T],\calU)$ and $x_0\in\dom A$ be arbitrary. Then for $t\in [0,T]$ the classical solution $x := x_u(\,\cdot\,;x_0)$ satisfies
\begin{align*}
\tfrac{d}{dt}H(x(t))
&= \tfrac 12\big(\<\dot x,x\> + \<x,\dot x\>\big) = \Re\<x,\dot x\>\\
&= \Re\<x,Ax+Bu\> = \Re\<B^*x,u\> + \mathbf r[x]\\
&= \Re\<u,y\> - \|R^{1/2}x\|^2.
\end{align*}
Integrating this over $[0,T]$ shows that $(u,x)$ satisfies the energy balance
\begin{equation}\label{e:eb_classic}
\tfrac 12\|x(T)\|^2 - \tfrac 12\|x_0\|^2 = \int_0^T\Re\<u,y\>\,dt - \int_0^T\|R^{1/2}x\|^2\,dt.
\end{equation}
Now, let $u\in L^1(0,T;\calU)$ and $x_0\in\calH$ be arbitrary. Then there exist sequences $(u_n)\subset C^1([0,T],\calU)$ and $(z_n)\subset\dom A$ such that $u_n\to u$ in $L^1(0,T;\calU)$ and $z_n\to x_0$ in $\calH$. Set $x_n(t) := x_{u_n}(t;z_n)$ and $x(t) := x_u(t;x_0)$. Then
\begin{align*}
\|x_n&(t) - x(t)\|\\
&= \left\|T(t)(z_n-x_0) + \int_0^t T(t-s)B(u_n(s)-u(s))\,ds\right\|\\
&\le C\big(\|z_n-x_0\| + \|B\|\|u_n-u\|_1\big),
\end{align*}
where $C = \sup_{s\in [0,t_2]}\|T(s)\|$. In particular, $x_n\to x$ uniformly on $[0,t_2]$.

Set $u_{nm} := u_n-u_m$ and $x_{nm} := x_n-x_m$. Then $x_{nm} = x_{u_{nm}}(\,\cdot\,;z_n-z_m)$. Hence, $(u_{nm},x_{nm})$ satisfies \eqref{e:eb_classic} and we obtain $\int_{0}^{t_2}\|R^{1/2}x_{nm}(t)\|^2\,dt\to 0$ as $n,m\to\infty$. That is, $(R^{1/2}x_n)$ is a Cauchy sequence in $L^2(0,T;\calH)$ and thus converges to some $w\in L^2(0,T;\calH)$. In particular, there exists a subsequence $(x_n')$ of $(x_n)$ such that $R^{1/2}x_n'(t)\to w(t)$ for a.e.\ $t\in [0,t_2]$. From the closedness of $R^{1/2}Q$ we infer that $x(t)\in\dom R^{1/2}Q$ and $R^{1/2}x(t) = w(t)$ for these $t$ so that $R^{1/2}x = w\in L^2(0,t_2;\calH)$. It is now clear that $(u,x)$ satisfies \eqref{e:eb_classic} and \eqref{e:RR}.

``$\Sla$'' Assume now that there exists $R\ge 0$ in $\calH$ such that \eqref{e:RR} holds. We shall show that the form $\mathbf b[x,y] := \<R^{1/2}x,R^{1/2}y\>$ with $\dom(\mathbf b) = \dom R^{1/2}$ is a closed extension of $\mathbf r$, and hence $\mathbf r$ is closable. By polarization, since $\mathbf b$ is closed by definition and $\dom A\subset\dom R^{1/2} = \dom(\mathbf b)$, it suffices to show that $\mathbf r[x] = \mathbf b[x]$ for $x\in\dom A$. 
For this, let $x_0\in\dom A$ be arbitrary and consider the trivial control $u\equiv 0$. Then $x(t) = T(t)x_0$ is a classical solution to $\dot x = Ax + Bu = Ax$, $x(0) = x_0$, and is thus contained in $\dom A$ for all $t\ge 0$. Differentiating the energy balance \eqref{e:eb_classical} with respect to $T$ shows that
\begin{align*}
\mathbf r[x(t)]
&= -\Re\<Ax(t),x(t)\> = -\Re\<\dot x(t),x(t)\>\\
&= -H'(x(t))\dot x(t) = \|R^{1/2}x(t)\|^2
\end{align*}
for all $t\ge 0$. In particular, for $t=0$ this yields $\mathbf r[z] = \mathbf b[z]$.\hfill $\square$
\end{pf}

\begin{ex}\label{ex:first}
(a) The operator $A$ in $\calH = L^2(0,1)$ from Example \ref{ex:ex} is maximal dissipative, but there is no $R$ such that \eqref{e:RR} is satisfied. Indeed, we have $\mathbf r[x,y] = \frac 12x(0)\ol{y(0)}$ and thus $\mathbf r[x] = \tfrac 12|x(0)|^2$. Let $x_n(\omega) := (1-\omega)^n$, $n\in\N$. Then $\|x_n\|_{L^2} = (2n+1)^{-1/2}\to 0$, $\mathbf r[x_n-x_m] = 0$, but $\mathbf r[x_n] = \frac 12$. Hence, $\mathbf r$ is not closable.

\smallskip\noindent
(b) If $A$ can be written as $A = J-R$ with $J=-J^*$ and $R = R^*\ge 0$ (examples are the heat equation [where $J=0$] or a Timoshenko beam [where $R$ is bounded], see \cite{PhilScha21}), then $\mathbf r$ is obviously closable as $\mathbf r[x] = \<Rx,x\>$, $x\in\dom J\cap\dom R$.
\end{ex}

\smallskip
\begin{rem}
We may reduce the closability of the form $\mathbf r$ to the closability of an operator. Indeed, since $A$ is maximal dissipative, so is $A-I$, and $A-I$ is invertible. Hence, the operator
$$
Q := -\tfrac 12\big((A-I)^{-1} + (A-I)^{-*}\big)
$$
is well defined, bounded and non-negative self-adjoint. It is now easily seen that
\[
\big\|Q^{1/2}(A-I)x\big\|^2 = \|x\|^2 + \mathbf r[x],\qquad x\in\dom A.
\]
Hence, $\mathbf r$ is closable if and only if $Q^{1/2}(A-I)$ is closable.
\end{rem}

The following corollary shows that the form $\mathbf r$ is closed for a large class of maximal dissipative operators.

\medskip
\begin{cor}\label{c:analytic}
If $A$ generates an analytic semigroup, then there exists a non-negative self-adjoint operator $R$ in $\calH$ such that \eqref{e:RR} holds for all $x_0\in\calH$ and $u\in L^2_{\rm loc}(\R_0^+,\calU)$.
\end{cor}
\begin{pf}
Indeed, in this case, $T := -A$ is an m-sectorial operator in the sense of \cite{kato}. The form $\tf[x,y] := \<Tx,y\>$, $x,y\in\dom T$, is then sectorial with some angle $\theta\in (0,\frac\pi 2)$ and closable (see Thm.\ VI.1.27 in \cite{kato}). Define two symmetric forms $\mathbf h$ and $\mathbf k$ by
\begin{align*}
{\mathbf h}[x,y] &:= \tfrac 12(\tf[x,y] + \ol{\tf[y,x]})\\
{\mathbf k}[x,y] &:= \tfrac 1{2i}(\tf[x,y] - \ol{\tf[y,x]}),
\end{align*}
respectively. Then $\tf = {\mathbf h} + i{\mathbf k}$ and, since $\tf[x]$ is in the sector $S_\theta$, we have $|{\mathbf k}[x]|\le(\tan\theta){\mathbf h}[x]$ for all $x\in\dom(\tf)$. Hence, if $x_n\to 0$ and ${\mathbf h}[x_n-x_m]\to 0$, then also ${\mathbf k}[x_n-x_m]\to 0$ so that $\tf[x_n-x_m]\to 0$, which implies $\tf[x_n]\to 0$ because $\tf$ is closable. In particular, ${\mathbf h}[x_n] = \Re\tf[x_n]\to 0$. Thus, $\mathbf r = -{\mathbf h}$ is closable.\hfill $\square$
\end{pf}

We would like to point out that it was shown in \cite[Theorem 3.1]{ArliTret20} that there exist generators $A$ of analytic semigroups (even with compact resolvent) for which $\dom(\frac 12(A+A^*)) = \dom A\cap\dom A^* = \{0\}$. Hence, in this case, the self-adjoint operator $R$ from Corollary \ref{c:analytic} is not related with $\frac 12(A+A^*)$ at all.

\section{Conclusion}
We considered the dissipation rate of infinite-dimensional linear input-state-output port-Hamiltonian systems with distributed control. We proved that there always is a linear operator mapping initial state $x_0\in\calH$ and control $u\in L^2(\R_0^+,\calU)$ to a function $\calR(x_0,u)\in L^2_{\rm loc}(\R_0^+,\calH_A)$ such that the dissipation rate can be expressed as $\|\calR(x_0,u)(t)\|_A^2$. We furthermore proved that the dissipation rate can be written as $\|R^{1/2}x_u(t;x_0)\|^2$ with a non-negative self-adjoint operator $R$ if and only if the form $\mathbf r$ is closable. We provided an example where this is not the case and proved that it always holds if $A$ generates an analytic semigroup. Future research consists of extending the results to port-Hamiltonian system nodes as introduced in \cite{PhilReis23}.










\end{document}